\documentstyle{article}
\title{A note on Abelian varieties embedded in quadrics}
\author{Luis Fuentes Garc{\'{\i}}a \thanks{Supported by EAGER.}}
\date{}

\newtheorem{teo}{Theorem}[section]

\newtheorem{lemma}[teo]{Lemma}

\newtheorem{rem}[teo]{Remark}

\def\Te{{\cal O}}

\def\P{{\bf P}}

\def\ov{\overline}

\def\qed{\hspace{\fill}$\rule{2mm}{2mm}$}

\def\lrw{{\longrightarrow}}

\begin{document}

\maketitle

{\bf Abstract:} We show that if $A$ is a $d$-dimensional abelian
variety in a smooth quadric of dimension $2d$ then $d=1$ and $A$
is an elliptic curve of bidegree $(2,2)$ on a quadric. This
extends a result of Van de Ven which says that $A$
only can be embedded in $\P^{2d}$ when $d=1$ or $2$. 

\vspace{0.1cm}

\section{Introduction.}

Let $A$ be a $d$-dimensional abelian variety embedded in $\P^N$.
It is well known that $2d\leq N$. Moreover, in \cite{van} Van de
Ven proved that the equality holds only when $d=1$ or $2$.

It is a natural question to study the possibilities for $d$ when
the abelian variety $A$ is embedded in any other smooth
$2d$-dimensional variety $V$. In particular, here we study the
embedding in smooth quadrics. We obtain the following result:

\begin{teo}
If $A$ is a $d$-dimensional abelian variety in a smooth quadric of
dimension $2d$ then $d=1$ and $A$ is an elliptic curve of bidegree
$(2,2)$ on a quadric.
\end{teo}

We will use similar methods to Van de Ven's proof. The calculation
of the self intersection of $A$ in the quadric and the
Riemann-Roch theorem for abelian varieties allow only the cases
$d=1,2,3$.

The case $d=1$ is the classical elliptic curve of type $(2,2)$
contained in the smooth quadric of $\P^3$.

When $d=2$, $A$ is an abelian surface in $\P^5$. We see that is
the projection of an abelian surface $A'\subset \P^6$ given by a
$(1,7)$ polarization. By a result \cite{lazarsfeld} due to R.
Lazarsfeld, this is projectively normal and it is not contained in
quadrics. Therefore, $A$ is not contained in quadrics either.

Finally, two results \cite{iyer1}, \cite{iyer2} of J.N.Iyer allow
us to discard the case $d=3$.

\section{Proof of the Theorem.}\label{teorema}

Let $j:A\hookrightarrow Q$ be an embedding of a $d$-dimensional
abelian variety into a $2d$-dimensional smooth quadric, with
$d>1$. The Chow ring of the smooth quadric in codimension $d$ is
generated by cocycles $\alpha$ and $\beta$ with the relations
$\alpha^2=\beta^2=1$, $\alpha.\beta=0$. Thus, $A$ will be
equivalent to $a \alpha + b \beta$ and
$$
A.A=a^2+b^2.
 \eqno{(1)}
$$
On the other hand, by the self-intersection formula
(\cite{hartshorne}, pag 431) we have $A.A=j_*c_d(N_{A,Q})$. To
obtain $c_d(N_{A,Q})$, let us consider the normal bundle sequence:
$$
0\lrw T_{A}\lrw j^*T_{Q}\lrw N_{A,Q}\lrw 0
$$
Since the tangent bundle of an abelian variety is trivial, we see
that $c(N_{A,Q})=j^*(c(T_Q))$. We compute the class of the tangent
bundle of a quadric in the following lemma:

\begin{lemma}
Let $i:Q\hookrightarrow \P^{n+1}$ be an $n$-dimensional smooth
quadric in $\P^{n+1}$. Then
$$
c(T_{Q})=(1+\ov{H})^{n+2}(1+2\ov{H})^{-1}
$$
where $\ov{H}=i^*H$ and $H$ is a hyperplane in $\P^{n+1}$.
\end{lemma}
{\bf Proof:} We have an exact sequence:
$$
0\lrw T_{Q}\lrw i^*T_{\P^{n+1}}\lrw N_{Q,\P^{n+1}}\lrw 0
$$
Since $Q$ is a hypersurface $N_{Q,\P^{n+1}}\cong \Te_Q(Q)\cong
\Te_Q(2\ov{H})$ and the total class of the normal bundle is $
c(N_{Q,\P^{n+1}})=1+2\ov{H} $. On the other hand, it is well known
that $c(T_{\P^{n+1}})=(1+H)^{n+2}$. Now, from the splitting
principle the claim follows. \qed

Let us apply this lemma to the previous situation. We obtain
$$
\begin{array}{rl}
{c(N_{A,Q})}&{=(1+h)^{2d+2}(1+2h)^{-1}=}\\
{}&{\sum_{k=0}^{2d+2} \left(^{2d+2}_{\;\;\, k}\right) h^k\sum_{l=0}^{\infty} (-2h)^{-l}}\\
\end{array}
$$
where $h=j^*\ov{H}$. In particular, the top class is
$$
c_d=F_d h^d, \mbox{ with } F_d=\sum_{k=0}^{d}
\left(^{2d+2}_{\;\;\, k}\right) (-2)^{(d-k)}.
$$
Substituting this into the self-intersection formula, we have:
$$
A.A=F_dj_*(j^*\ov{H}^d)=F_d\ov{H}^d.j_*A=F_d
(a\alpha+b\beta).\ov{H}^d=F_d(a+b).
$$
Combining this expression with $(1)$ we obtain the following
relation
$$
a^2+b^2=F_d(a+b)
 \eqno{(2)}
$$
or equivalently,
$$
(a-\frac{F_d}{2})^2+(b-\frac{F_d}{2})^2=\frac{F_d^2}{2}.
$$
We are interested in bounding the degree of $A$, when $(a,b)$
satisfy this equation. Note that this is a circle of center
$(\frac{F_d}{2},\frac{F_d}{2})$ and radius $\frac{F_d}{\sqrt{2}}$.
Since $deg(A)=a+b$, it is clear that the maximal degree is reached
when $(a,b)=(F_d,F_d)$, that is,
$$
deg(A)\leq 2F_d.
$$
On the other hand, the abelian variety is embedded in $Q\subset
\P^{2d+1}$. When $d>2$, by Van de Ven's Theorem, it spans
$\P^{2d+1}$. Furthermore, by the Riemann-Roch theorem for abelian
varieties, we know that $h^0(\Te_A(h))=\frac{deg(A)}{d!}$. Thus,
we have the following inequality:
$$
deg(A)\geq 2(d+1)!
$$
Comparing the two bounds we see that a sufficient condition for
the non-existence of the embedding $j$ is $F_d<(d+1)!$. Now,
$$
F_d=\sum_{k=0}^{d} \left(^{2d+2}_{\;\;\, k}\right) (-2)^{(d-k)}
\leq \sum_{k=0}^{d} \left(^{2d+2}_{\;\;\, k}\right) (2)^{d} \leq
2^d2^{2d+1}=2^{3d+1}.
$$
We see that $(d+1)!>2^{3d+1}\geq F_d$ when $d=17$. A simple
inductive argument shows that this holds if $d\geq 17$.

If $d\leq 17$, using the exact value of $F_d$, we see that
$(d+1)!>F_d$ for any $d>3$.

We conclude that the unique possibilities are $d=2$ or $d=3$.

First, suppose that $A$ is an abelian surface contained in a
quadric. $F_2=7$ and we can check that the unique positive integer
solution of the equation $(2)$ is $a=b=7$. Thus $A$ must be an
abelian surface of degree $14$ given by the polarization $(1,7)$.
Note that $A\subset Q\subset \P^5$ is not linearly normal, that
is, it is the projection of a linearly normal abelian surface
$A'\subset \P^6$. The quadric $Q$ can be lifted to a quadric
containing the surface $A'$.

Lazarsfeld proved in \cite{lazarsfeld} that a very ample divisor
of type $(1,d)$ with $d\geq 13$ or $d=7,8,9$ is projectively
normal. From this the following sequence is exact:
$$
0\lrw H^0(I_{A',\P^6}(2))\lrw H^0(\Te_{\P^6}(2))\lrw
H^0(\Te_{A'}(2))\lrw 0
$$
Since $h^0(\Te_{\P^6}(2))=h^0(\Te_{A'}(2))=28$, there are not
quadrics containing the abelian surface $A'$ and we obtain a
contradiction.

Finally, suppose that $d=3$. Now, $F_3=24=(3+1)!$, so the degree
of the abelian variety is exactly $2F_3=48$. The line bundle
$\Te_{A}(h)$ corresponds to a divisor of type $(1,1,8)$ or
$(1,2,4)$. But J.N.Iyer prove in \cite{iyer1} that a line bundle
of type $(1,\ldots,1,2d+1)$ is never very ample. Moreover, in
\cite{iyer2} she studies the map defined by a line bundle of type
$(1,2,4)$ in a generic abelian threefold. She obtains that it is
birational but not an isomorphism onto its image. Note that the
very ampleness is an open condition for polarized abelian
varieties (see \cite{dehusp}). It follows that a linear system of
type $(1,2,4)$ cannot be very ample on any abelian threefold and
this completes the proof.

\begin{rem}
The sequence $F_d=\sum_{k=0}^{d} \left(^{2d+2}_{\;\;\, k}\right)
(-2)^{(d-k)}$ is related to the Fine numbers. For a reference see
\cite{deutsch}.
\end{rem}

{\bf Acknowledgement} I thank K. Hulek for suggesting me this
problem and for his advice and interest. I am grateful to the
Institut f\"{u}r Mathematik of Hannover for its hospitality and
especially to E. Schellhammer for his patience. I want also thank
E. Deutsch for his remark about bounding $F_d$ and its relation
with the Fine numbers.

{\bf E-mail:} luisfg@usc.es

Luis Fuentes Garc\'{\i}a. Institut f\"{u}r Mathematik. 30167
Hannover, Germany.

\end{document}